\documentclass[12pt,leqno]{article}
\usepackage{amsfonts}
\usepackage{amsmath, amsthm, amsfonts, amssymb, color}
\usepackage{mathrsfs}
\usepackage[notcite,notref]{showkeys}

\newtheorem{thm}{Theorem}[section]

\theoremstyle{definition}

\newcommand{\scr}[1]{\mathscr #1}
\definecolor{wco}{rgb}{0.5,0.2,0.3}

\numberwithin{equation}{section} \theoremstyle{remark}

\newcommand{\ua}{\uparrow}

\title{{\bf Talagrand Inequality on Free Path Space and Application to  Stochastic Reaction Diffusion Equations}\footnote{Supported in
 part by    NNSFC(11671372,  11431014, 11721101, 11771326, 11831014)}
}
\author{
{\bf   Feng-Yu Wang$^{a),c)}$,  Tusheng Zhang$^{b),d}$}\\
\footnotesize{$^{a)}$Center for Applied Mathematics, Tianjin University, Tianjin 300072, China}\\
\footnotesize{$^{b)}$School of Mathematics, University of Manchester, Oxford Road, Manchester
M13 9PL,   U.K.  }\\
 \footnotesize{$^{c)}$Department of Mathematics,
Swansea University, Singleton Park, SA2 8PP, UK}\\
\footnotesize{$^{d)}$ School of Mathematics,
University of Science and Technology of China, Hefei, China}\\
\footnotesize{wangfy@tju.edu.cn,  tusheng.zhang@manchester.ac.uk}}
\begin{document}
\def\R{\mathbb R}  \def\ff{\frac} \def\ss{\sqrt} \def\B{\mathbf
B} \def\W{{\mathbb W}}
\def\N{\mathbb N} \def\kk{\kappa} \def\m{{\bf m}}
\def\dd{\delta} \def\DD{\Dd} \def\vv{\varepsilon} \def\rr{\rho}
\def\<{\langle} \def\>{\rangle} \def\GG{\Gamma} \def\gg{\gamma}
  \def\nn{\nabla} \def\pp{\partial} \def\EE{\scr E}
\def\d{\text{\rm{d}}} \def\bb{\beta} \def\aa{\alpha} \def\D{\scr D}
  \def\si{\sigma} \def\ess{\text{\rm{ess}}}
\def\beg{\begin} \def\beq{\begin{equation}}  \def\F{\scr F}
\def\Ric{\text{\rm{Ric}}} \def\Hess{\text{\rm{Hess}}}
\def\e{\text{\rm{e}}} \def\ua{\underline a} \def\OO{\Omega}  \def\oo{\omega}
 \def\tt{\tilde} \def\Ric{\text{\rm{Ric}}}
\def\cut{\text{\rm{cut}}} \def\P{\mathbb P} \def\ifn{I_n(f^{\bigotimes n})}
\def\C{\scr C}      \def\aaa{\mathbf{r}}     \def\r{r}
\def\gap{\text{\rm{gap}}} \def\prr{\pi_{{\bf m},\varrho}}  \def\r{\mathbf r}
\def\Z{\mathbb Z} \def\vrr{\varrho} \def\l{\lambda}
\def\L{\scr L}\def\Tt{\tt} \def\TT{\tt}\def\II{\mathbb I}
\def\i{{\rm in}}\def\Sect{{\rm Sect}}\def\E{\mathbb E} \def\H{\mathbb H}
\def\M{\scr M}\def\Q{\mathbb Q} \def\texto{\text{o}} \def\LL{\Lambda}
\def\Rank{{\rm Rank}} \def\B{\scr B} \def\i{{\rm i}} \def\HR{\hat{\R}^d}
\def\to{\rightarrow}\def\l{\ell}\def\ll{\lambda}
\def\8{\infty}\def\ee{\epsilon}

\maketitle

\begin{abstract}
By using a split argument due to \cite{BWY}, the transportation cost inequality is established on the free path space of Markov processes. The general result is applied to  stochastic reaction diffusion equations with random initial values.
  \end{abstract}


\section{Introduction}
Let $(E,\rr)$ be a  metric space, and let $\scr P(E)$ be the class of all probability measures on $E$.    The  quadratic Warsserstein distance between $\mu_1,\mu_2\in \scr P(E)$ is defined
 by \begin{equation*}
\W_{2}(\mu_1,\mu_2)= \inf_{\pi\in\C(\mu_1,\mu_2)}\bigg\{\int_{E\times
E}\rr^2(x,y)\pi(\d x,\d y)\bigg\}^{1/2},
\end{equation*}
where $\mathscr{C}(\mu_1,\mu_2)$ is the space of all couplings of
$\mu_1$ and $\mu_2$. In the study of Monge-Kontorovich optimal transportation problem, this distance is explained as the minimal cost to transport   distribution $\mu_1$ into $\mu_2$ at the cost rate   (cost function) $\rr$.
Thus, an inequality involving $\W_2$ is called a transportation cost inequality (TCI). Since the optimal transportation is usually unknown, in applications it is important  to estimate
 $\W_2$ by easier to calculate quantities, for instance the relative entropy $H(\mu_1|\mu_2):=  \int_E \big(\log \ff{\d\mu_1}{\d\mu_2}\big)\d\mu_1$ if $\mu_1$ is absolutely continuous with respect to $\mu_2$, and $H(\mu_1|\mu_2):=\infty$ otherwise.

 In 1996,   Talagrand  \cite{T96} established the following beautiful TCI    for the standard Gaussian measure   $\mu$ on $\R^d$  with  $\rr(x,y)=|x-y|$:
\begin{equation*}
\W_{2}(\nu, \mu)^2\le 2 H(\nu|\mu),\ \ \nu\in \scr P(\R^d),
\end{equation*}
where the constant $2$ is sharp.   Since then, this type  TCI  has been intensively investigated and applied for
various different distributions, and   was linked to   functional inequalities, concentration
phenomena, optimal transport problem, and large deviations, see
\cite{BGL, BG99, GRS11, FS, L, OV, U, W04b} and references
therein. Moreover,   Talagrand type TCI    has also been established  on the
path spaces of stochastic processes, see e.g. \cite{DGW, WZ04,
WZ06} for diffusion processes on $\R^d$,  \cite{P11} for
multidimensional semi-martingales, \cite{BWY, U10}  for   stochastic differential equations (SDEs) with memory, \cite{FS, FW, W04, W04b, W11} for
(reflecting) diffusion processes on Riemannian manifolds,  \cite{W10}   for SDEs
driven by pure jump processes,    and  \cite{M10, S11}
for SDEs with L\'evy   or fractional noises.

Recently,  by using the Girsanov transformation argument developed from \cite{DGW}, the Talagrand inequality was established
on the path space for solutions of stochastic reaction diffusion equations with deterministic initial values, see \cite{KS}, \cite{SZ}. In this paper, we aim to extend this result to the case with random initial values. In this case, the distribution of a solution is a probability measure on the free path space, where the initial value is not fixed.
Since the Girsanov transformation does not change initial distributions, it  does not work  for   probability measures  with different initial distributions. However, two equivalent probability measures on the free path space may have different initial distributions.   To overcome this difficulty,    we will adopt a split argument used in \cite{BWY} to reduce the problem to the case with deterministic initial value, to which the Girsanov transformation applies.

The remainder of the paper is organized as follows. In Section 2 we present a general result   on the TCI for Markov processes with random initial values, which is then applied  in Section 3 to stochastic reaction diffusion equations.

   \section{A general result }
Let $(E,\rr)$ be a Polish space, and let $(P_t)_{t\ge 0}$ be the semigroup of a continuous Markov process on $E$.
For any $T>0$ and $\mu\in \scr P(E)$, let $P^\mu$ denote the distribution   of the Markov process  up to time $T$ with initial distribution $\mu$; i.e. letting $P_t(x,\cdot)$ be the associated Markov transition kernel,
$P^\mu$ is the unique probability measure on the free path space
$$ E_T:= C([0,T];E)\ \text{equipped\ with}\ \rr_T(\xi,\eta):= \sup_{t\in [0,T]}\rr(\xi_t,\eta_t),$$
such that for any $0=t_0<t_1\cdots<t_n=T$ and $\{A_i\}_{0\le i\le n}\subset \B(E)$,
$$P^\mu(X_{t_i}\in A_i, 0\le i\le n)= \int_{A_0} \mu(\d x_0)\int_{A_1}P_{t_1-t_0}(x_0,\d x_1)\cdots\int_{A_n} P_{t_n-t_{n-1}}(x_{n-1},\d x_n),$$
where $X_t, t\geq 0$ denotes the canonical  coordinate process on the path space $E_T$.
When $\mu=\dd_x$, the Dirac measure at $x\in E$, we simply denote $P^\mu=P^x.$ Then
\beq\label{MK} P^\mu= \int_E P^x\mu(\d x),\ \ \mu\in \scr P(E).\end{equation}

Let $\W_2$ and $\W_{2,T}$ be the Wasserstein distances induced by $\rr$ on $\scr P(E)$  and $\rr_T$ on $\scr P(E_T)$ respectively.
We aim to establish the TCI for $P^\mu$ by using those for $\{P^x: x\in E\}$ and $\mu$.

\beg{thm}\label{T1} Assume that for some   constants $c_1,c_2\in (0,\infty)$ one has
\beq\label{TM2} \W_{2,T}(Q,P^x)^2\le c_1 H(Q|P^x),\ \ x\in E, Q\in \scr P(E_T),\end{equation}
\beq\label{TM3} \W_{2,T}(P^x,P^y)^2\le c_2 \rr(x,y)^2,\ \ x,y\in E. \end{equation}
If $\mu\in \scr P(E)$ satisfies
\beq\label{TM1}   \W_{2}(\nu,\mu)^2\le c_0 H(\nu|\mu),\ \ \nu\in \scr P(E) \end{equation} for some constant $c_0\in (0,\infty)$, then
\beq\label{TM0} \W_{2,T}(Q,P^\mu)^2\le C H(Q|P^\mu),\ \ Q\in \scr P(E_T) \end{equation}
holds for $C= \big(\ss{c_1}+\ss{c_0c_2}\big)^2.$ On the other hand, $\eqref{TM0}$ implies $\eqref{TM1}$ for $c_0=C$.
\end{thm}

\beg{proof} (1) We first deduce \eqref{TM0} from \eqref{TM1}. Let $Q=FP^\mu\in \scr P(E_T)$ and   $u_0: E_T\to E$ with $u_0(\xi)=\xi_0$.
Then
\beq\label{QQ}  \{Q\circ u_0^{-1}\}(\d x)= p(x)\mu(\d x)=:\nu(\d x)\end{equation} holds for
$$p(x):= \int_{E_T} F(\xi) P^x(\d \xi),\ \ x\in E.$$  By the triangle inequality,
\beq\label{TRA} \W_{2,T}(Q,P^\mu)\le \W_{2,T}(Q, P^\nu)+ \W_{2,T}(P^\nu, P^\mu).\end{equation}
Below we estimate these two terms   respectively.

To estimate  $\W_{2,T}(Q, P^\nu),$  we note that  \eqref{MK} implies
 \beg{align*} &\int_{E_T}f(\xi_0) F(\xi) P^\mu(\d \xi)=  \int_E f(x) \mu(\d x) \int_{E_T} F(\xi) P^x(\d \xi)\\
 &=\int_E f(x) p(x) \mu(\d x) =\int_{E_T} (fp)(\xi_0) P^\mu(\d\xi),\ \ f\in \B_b(E).\end{align*}
 Therefore, letting $\E^\mu$ be the expectation with respect to $P^\mu$, we have
 \beq\label{QQ2} p\circ u_0= \E^\mu(F|u_0).\end{equation}
Now, let
$$F_x(\xi)= 1_{\{p(x)>0\}}\ff{F(\xi)}{p(x)},\ \ x\in E, \xi\in E_T.$$ By \eqref{TM2}, if $p(x)>0$ then
$$\W_{2,T}(F_xP^x,P^x)^2\le c_1  P^x(F_x\log F_x).$$  So, for any $G,H\in \scr C,$ where
$$\C:= \big\{(G,H): G,H\in C_b(E_T), G(\xi)\le H(\eta)+\rr_T(\xi,\eta)^2\ \text{for}\  \xi,\eta\in E_T\big\},$$
 we have
 $$\int_{E_T} F_x G \d P^x-\int_{E_T} H\d P^x\le c_1\int_{E_T} (F_x\log F_x)\d P^x,\ \ p(x)>0.$$
 Integrating with respect to $\nu(\d x):= p(x)\mu(\d x)$ and using \eqref{MK}, we obtain
\beg{align*} &Q(G)-P^\nu(H)=\int_{E_T} G\d Q -\int_{E_T} H\d P^\nu \\
&= \int_E  \bigg\{\int_{E_T} F_x  G \d P^x-\int_{E_T} H\d P^x\bigg\} p(x) \mu(\d x)\\
&\le c_1\int_E\bigg\{\int_{E_T} (F_x\log F_x)\d P^x\bigg\}p(x) \mu(\d x)\\ &=c_1 \int_{E_T} \big\{F\log F- F\log \E^\mu(F|u_0)\big\}\d P^\mu\\
&= c_1 H(Q|P^\mu) - c_1 \E^\mu [F\log \E^\mu(F|u_0)] \le c_1 H(Q|P^\mu),\end{align*}where the last step is due to the fact that
\beg{align*} &\E^\mu [F\log \E^\mu(F|u_0)]= \E^\mu [\E^\mu(F|u_0)\log \E^\mu(F|u_0)]\\
&\ge\E^\mu [\E^\mu(F|u_0)]\log \E^\mu [\E^\mu(F|u_0)]= \E^\mu [F]\log  \E^\mu [F]= 0.\end{align*}
 Therefore, by the Kontorovich dual formula, we arrive at
 \beq\label{EE1} \W_{2,T}(Q,P^\nu)^2= \sup_{(G,H)\in \C} \big\{Q(G)-P^\nu(H)\big\}\le c_1 H(Q|P^\mu).\end{equation}

 On the other hand, by \eqref{TM3}, for any $(G,H)\in \C$ we have
 \beq\label{EE'} \int_{E_T} G\d P^x- \int_{E_T} H\d P^y\le c_2\rr(x,y)^2,\ \ x,y\in E.\end{equation}
 Let $\pi\in \C(\nu,\mu)$ be the optimal coupling such that
 $$\W_2(\nu,\mu)^2=\int_{E\times E} \rr(x,y)^2\pi(\d x,\d y).$$
 Integrating \eqref{EE'} with respect to $\pi(\d x,\d y)$, and applying \eqref{MK}, we obtain
 $$\int_{E_T} G\d P^\nu- \int_{E_T} H\d P^\mu=\int_{E\times E}\bigg\{\int_{E_T} G\d P^x- \int_{E_T} H\d P^y\bigg\}\pi(\d x,\d y)\le c_2 \W_2(\nu,\mu)^2.$$ Combining this with  the Kontorovich dual formula, and applying \eqref{TM1}, we arrive at
\beq\label{TM5} \W_{2,T}(P^\nu, P^\mu)^2\le c_2 \W_2(\nu,\mu)^2\le c_0c_2 \mu(p\log p).\end{equation}
 Since  \eqref{MK}, \eqref{QQ2} and Jensen's inequality imply
\beg{align*}&\mu(p\log p)=\int_{E_T} \big\{(p\circ u_0)\log p\circ u_0\big\}\d P^\mu \\
 &= \E^\mu [\E^\mu(F|u_0)\log \E^\mu(F|u_0)]\le \E^\mu[\E^\mu(F\log F|u_0)]= H(Q|P^\mu),\end{align*}
 it follows from \eqref{TM5} that
 $$\W_{2,T}(P^\nu, P^\mu)^2\le c_0c_2 H(Q|P^\mu).$$
 Combining this with \eqref{TRA} and \eqref{EE1}, we prove \eqref{TM0}

 (2) To deduce \eqref{TM1} from \eqref{TM0}, for $\nu=p\mu$ we take $Q= (p\circ u_0) P^\mu$. Let $\Pi\in \C(Q,P^\mu)$ be the optimal coupling such that
 $$\W_{2,T}(Q,P^\mu)^2=\int_{E_T\times E_T}\rr_T^2\,\d\Pi.$$ We have $\pi:=\Pi\circ(u_0,u_0)^{-1}\in \C(\nu,\mu)$, so that
 \beg{align*} &\W_2(\nu,\mu)^2\le \int_{E\times E} \rr^2\d\pi= \int_{E_T\times E_T} \rr^2(\xi_0,\eta_0)  \Pi(\d \xi,\d\eta)\\
 &\le \int_{E_T\times E_T} \rr^2_T(\xi,\eta)  \Pi(\d \xi,\d\eta)=\W_{2,T}(Q,P^\mu)^2.\end{align*}
 Combining this with \eqref{TM0} and noting that \eqref{MK} implies
 $$H(Q|P^\mu)= \int_{E_T}\big\{( p\circ u_0)\log p\circ u_0\big\}\d P^\mu= \int_E (p\log p)\d\mu=H(\nu|\mu),$$
 we derive \eqref{TM1} for $c_0=C.$

\end{proof}

\section{TCI for stochastic reaction diffusion equations with random initial values}
 Let $C_0([0,1])=\{u\in C([0,1]): u(0)=u(1)=0\}$.
 Consider the following SPDE on $C_0([0,1])$:
 \begin{align}\label{3.1}
\left\{
\begin{aligned}
 &\d u_t(x) =\frac{1}{2}u_t^{\prime\prime}(x)dt+b(u_t(x))dt+ \sigma(u_t(x))W(\d t,\d x), \quad x\in (0,1), \\
 &u_t\in C_0([0,1]), \quad\quad t\ge  0,
\end{aligned}
\right.
\end{align}
where   $W(\d t,\d x)$ is a space-time white noise on a complete  probability
space $(\Omega, \F, \P)$ with natural filtration $\F_t$ generated by
 the Brownian sheet $\{W(t,x): (t,x)\in [0, \infty)\times [0,1]\}$,    $u_0$ is a  $C_0([0,1])$-valued random variable independent of  $W$, and
$b, \sigma: \mathbb{R}\rightarrow \mathbb{R}$ are
locally bounded measurable functions. We say that an adapted,
continuous process  $\{u_t\}_{t\ge 0}$ on $C_0([0,1])$ is a solution to   (\ref{3.1}),  if $\P$-a.s.
\beq\label{3.2}\beg{split}
&\int_0^1u_t(x)\phi(x)\d x=\int_0^1u_0(x)\phi(x)\d x
+\frac{1}{2}\int_0^t\d s\int_0^1u_s(x)\phi^{\prime\prime}(x)\d x \\
&+\int_0^t\d s\int_0^1b(u_s(x))\phi(x)\d x+ \int_0^t\int_0^1\sigma(u_s(x))\phi(x)W(\d s,\d x),\quad t\ge 0,
\phi \in C_0^2([0, 1]),\end{split}\end{equation}
 where $C_0^2([0,1]):=\{\phi\in C^2([0,1]): \phi(0)=\phi(1)=0\}.$
According to \cite{W},   $u_t$ is a solution to  (\ref{3.1}) if and only if  $\P$-a.s.
\beq \label{3.3}\beg{split}
u_t(x)=&P_t u_0(x)+\int_0^tP_{t-s} \{b(u_s)\}(x) \d s
 + \int_0^t\int_0^1p_{t-s}(x,y)\sigma(u_s(y))W(\d s,\d y),\ \ t\ge 0,\end{split}
\end{equation}
where $P_t $ and $p_{t}(x,y)$ are the Dirichlet heat semigroup and   heat kernel generated by    $\frac{1}{2}\Delta $ on $[0, 1]$.

We will apply Theorem \ref{T1} to
$$E:= C_0([0,1]),\ \ E_T:= C([0,T]; E)=C([0,T]; C_0([0,1])),$$
and $P^\mu$ being the distribution of the solution $(u_t)_{t\in [0,T]}$ with initial distribution $\mu\in \scr P(E).$
 To this end, we    need the following assumption.
\begin{itemize}
\item[ {\bf (H)} ] $\si$ is bounded,   $b$ and $\si$ are Lipschitz continuous.
 \end{itemize}
According to \cite{W},    when $b$ and $\si$ are Lipschitz continuous,  (\ref{3.1}) admits a unique solution for any (random) initial value $u_0$ on $E$. The boundedness of $\si$ was used in \cite{SZ} to establish the TCI for solutions of \eqref{3.1} with deterministic initial values.

\begin{thm}\label{T2} Assume   {\bf (H)} and let   $\mu\in \scr P(E).$ Then
 \begin{equation}\label{TMM}
W_2(Q, P^{\mu})\leq C H(Q|P^{\mu}),\ \ Q\in \scr P(E_T)
\end{equation} holds for some constant $C>0$ if and only if
\beq\label{TMM2}
W_2(\nu, \mu)\leq c H(\nu|\mu),\ \ \nu\in \scr P(E)
\end{equation} holds for some constant $c>0$.
\end{thm}

\beg{proof} In the present case, we have
\beg{align*} &\rr(f,g)= \sup_{x\in [0,1]}|f(x)-g(x)|,\ \  f,g\in E:=C_0([0,1]),\\
&\rr_T(\xi,\eta)= \sup_{(t,x)\in [0,T]\times [0,1]}|\xi_t(x)-\eta_t(x)|,\ \ \xi,\eta\in E_T:=C([0,T]];E).\end{align*}
According to \cite{SZ}, \eqref{TM2} holds for some constant $c_1>0$. So, by Theorem \ref{T1}, it suffices  to verify \eqref{TM3}. Letting  $u_t^f$ be the unique solution of \eqref{3.1} with $u_0=f\in E:=C_0([0,1]),$  we only need to prove
 \begin{equation}\label{5.1}
\E\left[\sup_{(t,x)\in [0,T]\times [0,1]}|u^f_t(x)-u^g_t(x)|^2\right]\leq c_2\sup_{x\in [0,1]}|f(x)-g(x)|^2,\ \ f,g\in C_0([0,1])
\end{equation}
for some constant $c_2>0$. Indeed,  since the law of $(u_t^f,u_t^g)_{t\in [0,T]}$ is a coupling of $P^f$ and $P^g$,
we have
$$\W_{2,T}(P^f,P^g)^2\le \E [\rr_T(u^f,u^g)^2]= \E\left[\sup_{(t,x)\in [0,T]\times [0,1]}|u^f_t(x)-u^g_t(x)|^2\right].$$ Below we prove the estimate \eqref{5.1}.

By \eqref{3.3} we have \begin{align}\label{add 0302.1}
  \E \left[\sup_{(t,x)\in[0,T]\times[0,1]}|u_t^f(x)-u_t^g(x)|^2\right] \leq 3\rr(f,g)^2+ 3(I_1 + I_2) ,
\end{align}
where
\begin{align*}
 & I_1 :=  \E\left[\sup_{(t,x)\in[0,T]\times[0,1]}\left|\int_0^t\int_0^1 p_{t-s}(x,y)\big[b(u_s^f(y))-b(u_s^g(y))\big]\,\d s\d y\right|^2\right],\\
  &I_2 :=   \E\left[\sup_{(t,x)\in[0,T]\times[0,1]}\left|\int_0^t\int_0^1 p_{t-s}(x,y)\big[\sigma(u_s^f(y))-\sigma(u_s^g(y))\big] W (\d s,\d y)\right|^2\right].
\end{align*}
Noting that the Dirichlet heat kernel satisfies
$$\sup_{x\in [0,1]}\int_0^t\d s\int_0^1 p_{t-s}(x,y)^2 \d y\le \ff{\ss{2t}}{\ss \pi},\ \ t>0, $$
and due to {\bf (H)} we have
\beq\label{LL} |b(x)-b(y)|\lor |\si(x)-\si(y)|\le K|x-y|,\ \ x,y\in [0,1]\end{equation}
for some constant $K>0$,
by H\"o{}lder's inequality we obtain
 \beq\label{I}\beg{split}
 & I_1 \leq K^2 \E\Bigg\{\sup_{(t,x)\in[0,T]\times[0,1]}\bigg[\left(\int_0^t\int_0^1 p_{t-s}(x,y)^2\,\d s\d y\right)\\
 &\qquad\qquad \qquad\qquad  \times\left(\int_0^t\int_0^1 |u_s^f(y)-u_s^g(y)|^2\,dsdy\right)\bigg]\Bigg\}\\
  \leq & \sqrt{\frac{2T}{\pi}}K^2 \int_0^T  \E\left[\sup_{(r,y)\in[0,s] \times[0,1]}|u_r^f(y)-u^g_r(y)|^2\right]\,\d s.
\end{split}\end{equation}
To estimate  the term $I_2$, we recall the following inequality due to \cite{SZ}: for any $T,\vv>0$, there exists a constant $C_{T,\vv}>0$ such that for any adapted random field $\gg(t,x)$ with $\E[\sup_{(s,x)\in [0,t]\times [0,1]}|\gg(s,x)|^2]<\infty, t\ge 0$, we have
\beq\label{4.4-1}\beg{split}
& \E\left[\sup_{(s,x)\in [0,t]\times [0,1]}\bigg|\int_0^s\int_0^1p_{s-r}(x,y)\gg(r,y)W(\d r,\d y)\bigg|^2\right ] \\
&\leq  \varepsilon \E\Big[\sup_{(s,x)\in [0,t]\times [0,1]}|\gg(s,x)|^2\Big]+C_{T,\vv}\int_0^t\E\Big[\sup_{(r,x)\in [0,s]\times [0,1]}|\gg(r,x)|^2\Big]\d r,\ \ t\in [0,T].
\end{split}\end{equation}
 Applying this to $\gg(s,x)=\si(u_s^f(x))-\si(u_s^g(x))$ and using \eqref{LL}, we obtain
 that for any $\epsilon>0$,
\beq\label{term II}\beg{split}
  I_2 \leq & \epsilon \E\left[\sup_{(t,x)\in[0,T]\times[0,1]}|\sigma(u_t^f(x))-\sigma(u_t^g(x))|^2\right] \\
  & + C_{T,\vv} \E\int_0^T\sup_{y\in[0,1]}\left|\sigma(u_s^f(y))-\sigma(u_s^g(y))\right|^2\,\d s  \\
  \leq & \epsilon K^2 \E\left[\sup_{(t,x)\in[0,T]\times[0,1]}|u^f_t(x)-u^g_t(x)|^2\right]  \\
  & + C_{T,\vv}K^2 \int_0^T \E\left[\sup_{(r,y)\in[0,s]\times[0,1]}\left|u^f_r(y)-u^g_r(y)\right|^2\right]\,\d s, \ t\in [0,T].
\end{split}\end{equation}
So, setting
$$
  Y(t):= \E\left[\sup_{(s,x)\in[0,t]\times[0,1]}|u_s^f(x)-u^g_s(x)|^2\right],
$$ which is finite for all $t\in [0,\infty)$ due to assumption {\bf (H)},
by combining  (\ref{add 0302.1})-(\ref{term II}) together  we obtain
$$
  Y(t) \leq  3\rr(f,g)^2+  3\sqrt{\frac{2T}{\pi}}K^2 \int_0^t Y(s)\,\d s + 3 \epsilon K^2 Y(t) +3 C_{T,\epsilon}K^2 \int_0^t Y(s)\,\d s,\ \ t\in [0,T].
$$
Choosing   $\varepsilon=\ff 1 {6K^2}$, we find a constant $c(T)>0$ such that
 $$Y(t) \leq 6\rr(f,g)^2+c(T) \int_0^t Y(s)\,\d s,\ \ t\in [0,T].$$
By Gronwall's inequality and $Y(t)<\infty$ for $t\ge 0$, this implies   \eqref{5.1} for $c_2=6\e^{c(T)T}.$

\end{proof}

  To illustrate Theorem \ref{T2}, we present examples of $\mu$ satisfying \eqref{TMM2}, such that
  \eqref{TMM} holds true.
  By \cite[Theorem 3.1]{FS},   the heat   measure on the loop space $C_0([0,1])$  satisfies \eqref{TMM2}.
  Next, by  Gross \cite{Gross}, the log-Sobolev inequality holds for the Brownian bridge measure $\mu_0$ on $C_0([0,1])$:
  $$\mu_0(F^2\log F^2)\le 2T \mu_0(\|DF\|_H^2)^2, \ \ F\in \D(D),\mu_0(F^2)=1,$$
  where $(D,\D(D))$ is the Malliavin gradient operator and $\|h\|_H:=(\int_0^T |h'_t|^2\d t)^{\ff 1 2}$ is the Cameron-Martin norm.
 So, by a standard perturbation argument, the log-Sobolev inequality
 $$\mu(F^2\log F^2)\le 2T \e^{{\rm osc}(V)} \mu(\|DF\|_H^2)^2, \ \ F\in \D(D),\mu(F^2)=1,$$
 holds for any probability measure $\d\mu=\e^{V}\d\mu_0$ with $V\in \B_b(C_0([0,1]))$, where ${\rm osc}(V):=\sup V-\inf V$. According to \cite[Theorem 1.10]{Shao}, this implies
  $$\tt\W_2(\nu,\mu)^2\le 2T \e^{{\rm osc}(V)} H(\nu,\mu),\ \ \nu\in \scr P(C_0([0,1])),$$
 where $\tt\W_2$ is the Wasserstein distance induced by   the   Cameron-Martin distance on $E$. Since   the Cameron-Martin distance is larger than the uniform distance $\rr$, \eqref{TMM2} holds for this class of measures $\mu$.

\end{document}